\documentclass[12pt]{amsart}
\usepackage{amssymb,upref,enumerate}

\usepackage{amsmath,amssymb,amscd,amsthm,amsxtra}
\usepackage{latexsym}
\usepackage{graphics,epsfig}
\allowdisplaybreaks

\def\XXint#1#2#3{{\setbox0=\hbox{$#1{#2#3}{\int}$}
\vcenter{\hbox{$#2#3$}}\kern-.5\wd0}}

\newtheorem{theorem}{Theorem}
\newtheorem{proposition}{Proposition}
\newtheorem*{theorem*}{Theorem}

\newtheorem{lemma}{Lemma}

\newtheorem{corollary}{Corollary}

\theoremstyle{remark}

\def\({\left(}
\def\){\right)}
\def\be {\begin{equation}}
\def\en{\end{equation}}

\def\Cdot{{\dot{C}}}
\def\div{\text{div}~}

\newcommand{\tensor}{\otimes}

\numberwithin{equation}{section}

\begin{document}

%\dedicatory{}



%\thanks{}

\title[LANS equation with rough data]{Lagrangian Averaged Navier-Stokes equation with rough data in Sobolev spaces}
\author{Nathan Pennington}
\address{Nathan Pennington, Department of Mathematics, Kansas State University, 138 Cardwell Hall,
Manhattan, KS-66506, USA.}

\email{npenning@math.ksu.edu}
\date{\today}

\begin{abstract}
The Lagrangian Averaged Navier-Stokes equation is a recently derived approximation to the Navier-Stokes equation.  In this article we prove the existence of short time solutions to the incompressible, isotropic Lagrangian Averaged Navier-Stokes equation with low regularity initial data in Sobolev spaces.  For $L^2$-based Sobolev spaces, we obtain global existence results.  More specifically, we achieve local existence with initial data in the Sobolev space $H^{n/2p,p}(\mathbb{R}^n)$.  For initial data in $H^{3/4,2}(\mathbb{R}^3)$, we obtain global existence, improving on previous global existence results, which required data in $H^{3,2}(\mathbb{R}^3)$.
\end{abstract}

\maketitle

\bigskip

\section{Introduction}
The Navier-Stokes equation describes the motion of many fluids, including ocean currents, the weather and water flowing through a pipe.  The Navier-Stokes equation is given by
\begin{equation}\label{NS}\aligned \partial_t u + (u\cdot \nabla)u=\nu\triangle u-\nabla p,
\\ u(0,x)=u_0(x),
\endaligned
\end{equation}
where $u:I\times \mathbb{R}^n\rightarrow \mathbb{R}^n$ for some time strip $I=[0,T)$, $\nu>0$ is a constant due to the viscosity of the fluid, $p:I\times\mathbb{R}^n\rightarrow \mathbb{R}^n$ denotes the fluid pressure, and $u_0:\mathbb{R}^n\rightarrow\mathbb{R}^n$.  The requisite differential operators are defined by $\triangle=\sum_{i=1}^n \frac{\partial^2}{\partial_{x_i}^2}$ and $\nabla=\(\frac{\partial} {\partial_{x_i}},...,\frac{\partial}{\partial_{x_n}}\)$.  In this article we consider incompressible fluids, which imposes the additional requirement that $\div u=\div u_0=0$, where $\div=\sum_{i=1}^n \frac{\partial}{\partial_{x_i}}$.

In dimension $n=2$, local and global existence of solutions to the Navier-Stokes equation are well known (see \cite{lady356}; for a more modern reference, see Chapter $17$ of \cite{T3}).  For dimension $n\geq 3$, the problem is significantly more complicated.  There is a robust collection of local existence results, including \cite{Kato}, in which Kato proves the existence of local solutions to the Navier-Stokes equation with initial data in $L^n(\mathbb{R}^n)$; \cite{KP}, where Kato and Ponce solve the equation with initial data in the Sobolev space  $H^{n/p-1,p}(\mathbb{R}^n)$; and \cite{KT}, where Koch and Tataru establish local existence with initial data in the space $BMO^{-1}(\mathbb{R}^n)$ (for a more complete accounting of local existence theory for the Navier-Stokes equation, see \cite{nsbook}).  In all of these local results, if the initial datum is assumed to be sufficiently small, then the local solution can be extended to a global solution.   However, the issue of global existence of solutions to the Navier-Stokes equation in dimension $n\geq 3$ for arbitrary initial data remains a significant open problem.

Because of the intractability of the Navier-Stokes equation,
several different equations that approximate the Navier-Stokes
equation have been studied. A recently derived approximating
equation is the Lagrangian Averaged Navier-Stokes equation (LANS).
The LANS equation is derived from the Lagrangian Averaged Euler (LAE)
equation analogously to the way the Navier-Stokes equation is derived from
the Euler equation.  The Euler equation itself is derived by computing the critical points of the energy functional.  Similarly, the curves satisfying the LAE equation are critical points of an averaged version of the energy functional, with the averaging occurring at the level of the initial data.  For an exhaustive treatment of this process, see \cite{Shkoller}, \cite{SK}, \cite{MRS} and \cite{MS2}. In \cite{MKSM} and \cite{CHMZ}, the authors discuss the numerical improvements that use of the LANS equation provides over more common approximation techniques of
the Navier-Stokes equation.

Like the Navier-Stokes equation, the LANS equation has both a
compressible and an incompressible formulation.  The compressible
LANS equation is derived and studied in \cite{compM}.  The
incompressible LANS equation exists most generally in the
anisotropic form, and is derived and studied in \cite{MS2}.  In this paper, we will
consider a special case of the anisotropic equation called the
isotropic incompressible LANS equation.  One form of the
incompressible, isotropic LANS equation on a region without
boundary is
\begin{equation}\label{LANS}\aligned\partial_t u+(u\cdot\nabla)u+\div
\tau^\alpha u=-(1-\alpha^2\triangle)^{-1}\nabla p+\nu\triangle u
\\ u(0,x)=u_0,~~ \div u=\div u_0=0
\endaligned
\end{equation}
where, as in $(\ref{NS})$, $u$ and $p$ map $I\times \mathbb{R}^n$ into $\mathbb{R}^n$, $u_0:\mathbb{R}^n\rightarrow\mathbb{R}^n$, $\alpha>0$ is a constant resulting from the Lagrangian averaging, and the Reynolds stress $\tau^\alpha$ is given by
\begin{equation*}\tau^\alpha u=\alpha^2(1-\alpha^2\triangle)^{-1}[Def(u)\cdot
Rot(u)],
\end{equation*}
where $Rot(u)=(\nabla u-\nabla u^T)/2$ and $Def(u)=(\nabla u+\nabla u^T)/2$.

For the LANS equation, establishing local existence is more complicated than in the Navier-Stokes setting, due to the presence of the additional non-linear term.  Specifically, for the LANS equation, we must (essentially) estimate the norm of both $\nabla u \nabla u$ and $u\tensor u$ (see $(\ref{LANS2})$ for the notation), while only the second term must be addressed for the Navier-Stokes equation.  However, long time behavior of a known local solution to the LANS equation is easier to control than in the Navier-Stokes setting.  In \cite{MS}, the authors establish local and global existence of solutions to the LANS equation with arbitrary initial data in the Sobolev space $H^{3,2}(\mathbb{R}^3)$.

In this paper, we seek local solutions to the LANS equation with initial data Sobolev spaces with minimal regularity and, when possible, we extend these local solutions to global solutions.  In particular, we obtain global solutions with initial data in $H^{s,2}(\mathbb{R}^3)$, where $s\geq 3/4$, improving on the result from \cite{MS}.

The paper is organized as follows.  The main results in this article are Theorems $\ref{special case 1}$ through $\ref{globathm2}$, and their statements can be found in Section $\ref{Main Theorems}$ below.  Theorems $\ref{special case 1}$ through $\ref{thm2}$ are local existence results.  Theorem $\ref{globathm2}$ states an $\emph{a priori}$ estimate used to extend some of the local solutions to global solutions, and the statement of the global results is given Corollary $\ref{globathm}$.  The proofs of Theorems $\ref{special case 1}$ and $\ref{thm1}$ are in Section $\ref{Solutions in the class of weighted continuous functions in time}$, and the proofs of Theorems $\ref{special case 2}$ and $\ref{thm2}$ are in Section $\ref{Solutions in the class of integral norms in time}$.  The proof of Theorem $\ref{globathm2}$ can be found in Section $\ref{Global Existence in Sobolev space}$.  Section $\ref{Higher regularity for the local existence result}$ contains two technical propositions necessary for the application of Theorem $\ref{globathm2}$ to our local existence results.  

We conclude this section by establishing some necessary notation.  As mentioned above, Sobolev spaces will be denoted by $H^{s,p}(\mathbb{R}^n)$  with norm $\|\cdot\|_{s,p}$.  We define the space
\begin{equation*}C^T_{a;k,q}=\{f\in
C((0,T):H^{k,q}(\mathbb{R}^n)):\|f\|_{a;k,q}<\infty\},
\end{equation*}
where
\begin{equation*}\|f\|_{a;k,q}=\sup\{t^a\|f(t)\|_{k,q}:t\in
(0,T)\},
\end{equation*}
$T>0$, $a\geq 0$, and $C(A:B)$ is the space of continuous functions from $A$ to $B$.  We let $\Cdot^T_{a;k,q}$ denote the subspace of $C^T_{a;k,q}$ consisting
of $f$ such that
\begin{equation*}\lim_{t\rightarrow 0^+}t^a f(t)=0
~\text{(in}~H^{k,q}(\mathbb{R}^n)).
\end{equation*}
Note that while the norm $\|\cdot\|_{a;k,q}$ lacks an explicit reference to $T$, there is an implicit $T$ dependence.  We also say $u\in BC(A:B)$ if $u\in C(A:B)$ and $\sup_{a\in A}\|u(a)\|_{B}<\infty$.  Lastly, setting $\mathbb{M}((0,T):\mathbb{E})$ to be the set of measurable
functions defined on $(0,T)$ with values in the space $\mathbb{E}$, we define
\begin{equation*}L^a((0,T):H^{k,q}(\mathbb{R}^n))=
\{f\in \mathbb{M}((0,T):H^{k,q}(\mathbb{R}^n)):(\int_0^T
\|f(t)\|^a_{k,q}dt)^{1/a}<\infty\}.
\end{equation*}

\textbf{Acknowledgments:} This work has its origins in the author's graduate thesis, completed at the University of North Carolina-Chapel Hill.  The author would like to thank the UNC Mathematics Department, his dissertation committee, and especially his advisor Professor Michael Taylor for the helpful suggestions and support.

\section{Main Theorems}\label{Main Theorems}

In this section we state our main results, and we begin with our local existence theorems.  Note that, because we are considering low-regularity spaces, our solutions are necessarily weak solutions.  A more explicit accounting of this can be found in Section $\ref{Solutions in the class of weighted continuous functions in time}$.  Finally, all of these local solutions depend uniquely and continuously on the initial data.  

The first theorem addresses initial data in the space $\dot{C}^T_{a;k,q}$.  

\begin{theorem}\label{special case 1}Let $u_0\in H^{3/4,2}(\mathbb{R}^3)$ be divergence-free.  Then there exists a local solution $u$ to the LANS equation $(\ref{LANS})$, where
\begin{equation*}u\in BC([0,T):H^{3/4,2}(\mathbb{R}^3))\cap \dot{C}^T_{1/8;1,2},
\end{equation*}
and $T$ is a non-increasing function of $\|u_0\|_{3/4,2}$, with $T=\infty$ if $\|u_0\|_{3/4,2}$ is sufficiently small.

Next, for $p>n\geq 3$,let $u_0\in H^{n/(2p),p}(\mathbb{R}^n)$ be divergence-free.  Then there exists a local solution $u$ to the LANS equation $(\ref{LANS})$, where
\begin{equation*}u\in BC([0,T):H^{n/2p,p}(\mathbb{R}^n))\cap \dot{C}^T_{(1-n/2p)/2;1,p},
\end{equation*}
and $T$ is a non-increasing function of $\|u_0\|_{n/2p,p}$, with $T=\infty$ if $\|u_0\|_{n/2p,p}$ is sufficiently small.
\end{theorem}

Our second theorem is analogous to the first, with the space $L^a((0,T):H^{s,p}(\mathbb{R}^n))$ replacing $\dot{C}^T_{a;k,q}$.

\begin{theorem}\label{special case 2}Let $u_0\in H^{3/4,2}(\mathbb{R}^3)$ be divergence-free.  Then there exists a local solution $u$ to the LANS equation $(\ref{LANS})$, where
\begin{equation*}u\in BC([0,T):H^{3/4,2}(\mathbb{R}^3))\cap L^{8}((0,T):H^{1,2}(\mathbb{R}^n)),
\end{equation*}
and $T$ is a non-increasing function of $\|u_0\|_{3/4,2}$, with $T=\infty$ if $\|u_0\|_{3/4,2}$ is sufficiently small.

Next, for $p>n\geq 3$, let $u_0\in H^{\frac{n}{2p},p}(\mathbb{R}^n)$ be divergence-free.   Then there exists a local solution $u$ to the LANS equation $(\ref{LANS})$ where
\begin{equation*}u\in BC([0,T):H^{n/2p,p}(\mathbb{R}^n))\cap L^{2/(1-n/2p)}((0,T):H^{1,p}(\mathbb{R}^n)),
\end{equation*}
and $T$ is a non-increasing function of $\|u_0\|_{n/2p,p}$, with $T=\infty$ if $\|u_0\|_{n/2p,p}$ is sufficiently small.
\end{theorem}

We note that in each of these theorems, we achieve local existence for initial data with regularity $n/2p$.  In the next two theorems, the minimal regularity achieved is only $n/p$.  However, these two theorems are significantly more general, allowing for a much wider range of parameter values, specifically with regard to the auxiliary space.

\begin{theorem}\label{thm1}For any divergence-free $u_0\in H^{r,p}(\mathbb{R}^n)$, there exists a local solution $u$ to $(\ref{LANS})$ such that 
\begin{equation*}u\in BC([0,T):H^{r,p}(\mathbb{R}^n))\cap \dot{C}^T_{a;k,c},
\end{equation*}
for some $T>0$, provided there exists a real number $b'$ such that the list of
conditions $(\ref{conditions})$ is satisfied (where $r=n/p+b)$.  Also, $T=T(\|u_0\|_{n/2p,p})$ can be chosen to be a non-increasing function of $\|u_0\|_{r,p}$, with $T=\infty$ if $\|u_0\|_{r,p}$ is sufficiently small.  
\end{theorem}

\begin{theorem}\label{thm2}For any divergence-free $u_0=\in H^{r,p}(\mathbb{R}^n)$ there exists a local solution $u$ to $(\ref{LANS})$ such that
\begin{equation*}u\in BC([0,T):H^{r,p}(\mathbb{R}^n)) \cap L^a((0,T):H^{k,c}(\mathbb{R}^n)),
\end{equation*}
for some $T>0$, provided the parameters (with $r=n/p+b$) satisfy
$(\ref{conditions2})$.  Also, $T=T(\|u_0\|_{r,p})$ can be chosen to be a non-increasing function of $\|u_0\|_{r,p}$, with $T=\infty$ if $\|u_0\|_{r,p}$ is sufficiently small.
\end{theorem}

We remark that in \cite{KP}, the authors achieved local existence for the Navier-Stokes equation with initial data regularity $n/p-1$, and the source of the gap between their result and the results presented here is the additional non-linear term in the LANS equation.  

The final theorem is an $\emph{a~ priori}$ estimate, providing a uniform-in-time bound on solutions to the LANS equation.

\begin{theorem}\label{globathm2}Let $u$ be a solution to the LANS equation $(\ref{LANS})$ on $I\times \mathbb{R}^3$, where $I=[0,T)$.  Then
\begin{equation*}\|u(t,\cdot)\|_{2,2}\leq M,
\end{equation*}
where $M$ is a function $\|u_0\|_{1,2}$.
\end{theorem}
Combining any of the $L^2(\mathbb{R}^3)$-based local results from Theorems $\ref{special case 1}$ through $\ref{thm2}$ with Theorem $\ref{globathm2}$ and applying a standard extension argument proves the following corollary.

\begin{corollary}\label{globathm}Let $u_0\in H^{3/4,2}(\mathbb{R}^3)$ and let $u$ be the associated local solution given by Theorem $\ref{special case 1}$ or Theorem $\ref{special case 2}$.  Then the local solution $u$ can be extended to a global solution.  Similarly, for $u_0\in H^{3/2,2}(\mathbb{R}^3)$,  and $u$ the associated local solution given by Theorem $\ref{thm1}$ or Theorem $\ref{thm2}$,  the local solution $u$ can be extended to a global solution.
\end{corollary}

We remark that to apply Theorem $\ref{globathm2}$ to any of the local results, we require Propositions $\ref{higher regularity theorem}$ and $\ref{higher regularity theorem2}$ from Section $\ref{Higher regularity for the local existence result}$.

\section{Solutions in $\dot{C}^T_{a;k,q}$}\label{Solutions in the class of weighted continuous functions
in time}

We begin by re-writing the LANS equation as 
\begin{equation}\label{LANS2}\aligned\partial_t u-Au+P^\alpha (\div\cdot(u\tensor
u)+\div\tau^\alpha u)=0,
\endaligned
\end{equation}
where we set $\nu=1$ in $(\ref{LANS})$, $A=P^\alpha\triangle$, $u\tensor u$ is the tensor with
$jk$-components $u_ju_k$ and $\text{div}\cdot(u\tensor u)$ is the
vector with $j$-component $\sum_k\partial_k(u_ju_k)$.  $P^\alpha$ is
the Stokes Projector, defined as
\begin{equation*}P^\alpha(w)=w-(1-\alpha^2\triangle)^{-1}\nabla f,
\end{equation*}
where $f$ is a solution of the Stokes problem: Given $w$, there is a
unique divergence-free $v$ and a unique (up to additive constants) function $f$ such
that
\begin{equation*}(1-\alpha^2\triangle)v+\nabla f=(1-\alpha^2\triangle)w.
\end{equation*}
For a more explicit treatment of the Stokes
Projector, see Theorem 4 of \cite{SK}.

Using Duhamel's principle, we write ($\ref{LANS2}$) as the
integral equation
\begin{equation}\label{intversion}u=\Gamma\varphi-G\cdot P^\alpha(\div(u\tensor u+\tau^\alpha(u))) 
\end{equation}
with
\begin{equation*}(\Gamma\varphi)(t)=e^{tA}\varphi,
\end{equation*}
where $A$ agrees with $\triangle$ when restricted to $P^\alpha
H^{r,p}$, and  
\begin{equation*}G\cdot g(t)=\int_0^t e^{(t-s)A}\cdot
g(s)ds.
\end{equation*}

Our plan is to construct a contraction mapping based on $(\ref{intversion})$.  Section $\ref{preliminary work}$  contains the operator estimates necessary for this contraction.  The proof of Theorem $\ref{thm1}$ is in Section $\ref{Proof of thm1}$, and the proof of Theorem $\ref{special case 1}$ is in Section $\ref{Improvements for special cases}$.

\subsection{Preliminary work}\label{preliminary work}
We begin by examining the Reynolds stress term.

\begin{lemma}\label{lemmaold}Let $r\in[1,\infty)$ and $1<p,q<\infty$, with $2/p-1/q<1$ and $0\leq n(2q-p)/(pq)\leq r-1$.  Then $\div\tau^\alpha:H^{r,p}(\mathbb{R}^n)\rightarrow H^{r,q}(\mathbb{R}^n)$, where we recall $\tau^\alpha(u)=\alpha^2(1-\alpha^2\triangle)^{-1}(Def(u)\cdot Rot(u)).$  Specifically, we have the estimate
\begin{equation*}\|\div\tau^\alpha(u)\|_{r,q}\leq C\|u\|^2_{r,p}
\end{equation*}
\end{lemma}

\begin{proof}  Recalling the definitions of $Def(u)$ and $Rot(u)$ and applying Proposition $1.1$ from \cite{TT} (which has its origins in \cite{kp2}, \cite{CW}, and \cite{kenigponcevega}) we get
\begin{equation*}\aligned \|\tau^\alpha(u)\|_{H^{r+1,q}(\mathbb{R}^n)}&\leq C\|[Def(u)\cdot Rot(u)]\|_{r-1,q}
\\ &\leq C\|\nabla u\|^2_{r-1,p}
\\ &\leq C\|u\|^2_{r,p}.
\endaligned
\end{equation*}
Since the divergence is a degree one differential operator, we get
\begin{equation*}\|\div \tau^\alpha(u)\|_{r,q}\leq \|\tau^\alpha(u)\|_{r+1,q}\leq C\|u\|^2_{r,p},
\end{equation*}
which proves the lemma.
\end{proof}
This immediately gives the following corollary.
\begin{corollary}\label{cortau}With the parameters as in Lemma $\ref{lemmaold}$, we have that $\div:\tau^\alpha:\dot{C}^T_{a;r,p}\rightarrow \dot{C}^T_{2a;r,q}$ with the estimate
\begin{equation*}\|\div\tau^\alpha(u)\|_{2a;r,q}\leq C\|u\|^2_{a;r,p}.
\end{equation*}
\end{corollary}

The purpose of this previous calculation is to note that $r\geq 1$ is forced by the presence of the $\tau^\alpha$ term.  For the Navier-Stokes equation, the $\tau^\alpha$ term is not present, and the corresponding result in that setting only requires $r\geq 0$.  This is the reason our local existence theory requires one more point of regularity than in the Navier-Stokes setting.

Our next task is to establish some properties of the operator $V^\alpha$ defined by
\begin{equation}\label{V}V^\alpha(u,v)=\div(u\tensor v)+\div\tau^\alpha(u,v))).
\end{equation}
Abusing notation, we will write $V^\alpha(u)=V^\alpha(u,u)$.  We also observe that $V^\alpha$ in linear in each of its arguments.  The following follows directly from Corollary $\ref{cortau}$.

\begin{proposition}\label{Vprop}Let $a\geq 0$, $b\geq 0$, $1<q,p<\infty$,
$0\leq n(2q-p)/pq\leq b$ and $2/p-1/q<1$.  Then
\begin{equation*}V^\alpha:\dot{C}^T_{a;b,p}\times\dot{C}^T_{a;b,p}\rightarrow \dot{C}^T_{2a;b-1,q}
\end{equation*}
with the estimate
\begin{equation}\label{Vest}\|V^\alpha(u,v)\|_{2a;b-1,q}\leq \|u\|_{a;b,p}\|v\|_{a;b,p}.
\end{equation}
\end{proposition}

Next, we observe that
\begin{equation*}V^\alpha(u)-V^\alpha(v)=-(V^\alpha(u,u-v)+V^\alpha(u-v,v)).
\end{equation*}
Using Proposition $\ref{Vprop}$, we have
\begin{equation*}\|V^\alpha(u,u-v)\|_{b-1,q}\leq \|u\|_{b,p}\|u-v\|_{b,p}
\end{equation*}
and
\begin{equation*}\|V^\alpha(u-v,v)\|_{b-1,q}\leq \|v\|_{b,p}\|u-v\|_{b,p}.
\end{equation*}

These estimates give that
\begin{equation}\label{prelip}\|V^\alpha(u(s))-V^\alpha(v(s))\|_{b-1,q}\leq C\(\|u(s)\|_{b,p}+\|v(s)\|_{b,p}\)\|u(s)-v(s)\|_{b,p}.
\end{equation}
Multiplying both sides by $t^a$ and distributing through the right
hand side, we get the following corollary to Proposition
$\ref{Vprop}$.
\begin{corollary}\label{Vcor}With the same assumptions on the parameters as in Proposition
$\ref{Vprop}$, we have that if $u,v\in \dot{C}^T_{a/2;b,q}$ then
\begin{equation}\label{Vcont}
\|V^\alpha(u(s))-V^\alpha(v(s))\|_{a;b-1,q}\leq
C(\|u\|_{a/2;b,p}+\|v\|_{a/2;b,p})\|u-v\|_{a/2;b,p}.
\end{equation}
\end{corollary}

Our next topic is the operator $\Gamma$.
\begin{proposition}\label{gamma prop}Let $s'\leq s''$, $1<q'\leq q''<\infty$,
and define $k''=(n/q'-n/q''+s''-s')/2$. Then
\begin{equation*}\|\Gamma f\|_{k'';s'',q''}\leq C\|f\|_{s',q'}
\end{equation*}
and for any $\varepsilon>0$, there exists sufficiently small $T$ such that
\begin{equation*}\|\Gamma f\|_{k'';s'',q''}\leq \varepsilon,
\end{equation*}
provided $k''>0$.
\end{proposition}
This is an immediate consequence of equation $1.15$ in Chapter $15$ of \cite{T3}.

We now turn our attention to the operator $G$. Assuming $s'\leq
s''$, $q'\leq q''$, and $u\in\dot{C}^T_{k';s',q'}$, we formally
calculate
\begin{equation*}\aligned
\|G\cdot u\|_{s'',q''}&=\|\int_0^t e^{(t-s)A}u(s)ds\|_{s'',q''}
\\& \leq C\int_0^t\|e^{(t-s)A}u(s)\|_{s'',q''}ds
\\& \leq C\int_0^t (t-s)^{-(s''-s'+n/q'-n/q'')/2}\|u(s)\|_{s',q'}ds
\\& \leq C\int_0^t (t-s)^{z}s^{-k'}s^{k'}\|u(s)\|_{s',q'} ds
\\& \leq Ct^{z-k'+1}\|u\|_{k';s',q'}
\endaligned
\end{equation*}
where $z=-(s''-s'+n/q'-n/q'')/2.$  This result will hold provided
$0\leq(s''-s'+n/q'-n/q'')/2<1$ and $k'<1$, and this leads to our
first result involving $G$.

\begin{proposition}\label{G}With $s'\leq s''$, $q'\leq q''$ and setting
$k''=k'-1+(s''-s'+n/q'-n/q'')/2$, $G$ continuously maps
$\dot{C}^T_{k';s',q'}$ into $\dot{C}^T_{k'';s'',q''}$ with $0\leq
(s''-s'+m/q'-m/q'')/2<1$ and $k'<1$ with the estimate
\begin{equation*}\|G\cdot u\|_{k'';s'',q''}\leq
C\|u\|_{k';s',q'}.
\end{equation*}
\end{proposition}

\subsection{Proof of Theorem $\ref{thm1}$}\label{Proof of thm1} To
prove Theorem $\ref{thm1}$ we begin by constructing the nonlinear
map
\begin{equation*}\Phi u=\Gamma \varphi-G\cdot P^\alpha(\div(u\tensor u)+\div
\tau^\alpha u).
\end{equation*}
Our goal is to show that this map is a contraction on an appropriate
function space.  Using $(\ref{V})$, $\Phi$ can be re-written as
\begin{equation*}\Phi u=\Gamma \varphi-G\cdot P^\alpha(V^\alpha(u)).
\end{equation*}

Beginning with initial data $\varphi\in H^{r,p}(\mathbb{R}^n)$ where
$r=\frac{n}{p}+b$, we construct the space
\begin{equation*}E_{T,M}=\{v\in BC([0,T):H^{r,p}(\mathbb{R}^n)) \cap \dot{C}^T_{a;k,c}
:\|v-\Gamma\varphi\|_{0;r,p}+\|v\|_{a;k,c}\leq M\}.
\end{equation*}
Our goal will be to
show that $\Phi$ is a contraction on this space for appropriate
choices of parameters.

To show $\Phi$ is a contraction, we use the mapping
properties of $G$ to send each component space of $E_{T,M}$ into an
intermediate space, and then use Proposition $\ref{Vprop}$.  Our
intermediate space will be of the form $\dot{C}^T_{2a;k-b',\bar{c}}$.

Our first task is to show that $\Phi$ maps $E_{T,M}$ into $E_{T,M}$.
To do this, we need to estimate
\begin{equation}\label{Phi1}\|\Phi
(u)-\Gamma\varphi\|_{0;r,p}=\|G\cdot P^\alpha V^\alpha(u)\|_{0;r,p}
\end{equation}
and
\begin{equation}\label{Phi2}\|\Phi(u)\|_{a;k,c}=
\|\Gamma\varphi-G\cdot P^\alpha V^\alpha (u)\|_{a;k,c}.
\end{equation}

To estimate $(\ref{Phi1})$, we note that by Proposition $\ref{G}$
and that $P^\alpha$ is a projection, we have that
\begin{equation}\label{1}\|G V^\alpha u\|_{0;\frac{n}{p}+b,p}
\leq C\|V^\alpha u\|_{2a;k-b',\bar{c}}
\end{equation}
will hold provided
\begin{equation}\label{A}\aligned 0&=2a-1+
\(\frac{n}{p}+b-(k-b')+\frac{n}{\bar{c}}-\frac{n}{p}\)/2
\\ 2a&<1
\\ 0&\leq (n/p+b-(k-b')+n/\bar{c}-n/p)/2<1
\\ k-b'&\leq n/p+b
\\ \bar{c}&\leq p.
\endaligned
\end{equation}

Proposition $\ref{Vprop}$ gives
\begin{equation*}\|V^\alpha u\|_{2a;k-b',\bar{c}}\leq C\|u\|^2_{a;k,c}
\end{equation*}
provided
\begin{equation}\label{AA}\aligned k\geq 0
\\ b'\geq 1,
\\ c>1,
\\ \bar{c}=\frac{nc}{2n-s'c}
\\ 0\leq s'\leq k-1
\\ s'c<n.
\endaligned
\end{equation}

These combine to give our estimate on $(\ref{Phi1})$.  To estimate
$\|G\cdot P^\alpha V^\alpha (u)\|_{a;k,c}$, we have
\begin{equation}\label{2}\|G V^\alpha u\|_{a;k,c}\leq C\|V^\alpha
u\|_{2a;k-b',\bar{c}}\leq C\|u\|^2_{a;k,c}
\end{equation}
will hold provided
\begin{equation}\label{B}\aligned a&=2a-1+\(k-(k-b')+\frac{n}{\bar{c}}-\frac{n}{c}\)/2
\\2a&<1
\\ \bar{c}&\leq c
\\ 0&\leq (k-(k-b')+n/\bar{c}-n/c)/2<1.
\endaligned
\end{equation}

Using $(\ref{1})$ and $(\ref{2})$, we have
\begin{equation*}\|\Phi(u)-\Gamma\varphi\|_{0;r,p}+\|\Phi(u)\|_{a;k,c}\leq
C\|u\|^2_{a;k,c}+\|\Gamma\varphi\|_{a;k,c}.
\end{equation*}

By assumption, $u\in E_{T,M}$, so $\|u\|^2_{a;k,c}\leq M^2$.  So our
last task is to estimate $\|\Gamma\varphi\|_{a;k,c}$.  From
Proposition $\ref{gamma prop}$, we have that
\begin{equation}\label{gamma prop 2*}\Gamma:\dot{C}^T_{0;\frac{n}{p}+b,p}\rightarrow
\dot{C}^T_{a;k,c}
\end{equation}
if $a>0$, $k\geq \frac{n}{p}+b$, $c\leq p$, and
\begin{equation*}a=\(\frac{n}{p}-\frac{n}{c}+k-\(\frac{n}{p}+b\)\)/2
\end{equation*}
which simplifies to
\begin{equation}\label{gammainequality}2a=k-\frac{n}{c}-b.
\end{equation}
Because $\Gamma\varphi\in \dot{C}^T_{a;k,c}$, there exists a $T$,
depending only on $M$ and the norm of the initial data $\varphi$,
such that $\|\Gamma\varphi\|_{a;k,c}\leq M/2$. So by choosing a
sufficiently small $M$ and an appropriate $T$, we have that
$\Phi:E_{T,M}\rightarrow E_{T,M}$.

Now we seek to show that $\Phi$ is a contraction map.  Let $u,v\in
E^{T,M}$.  Then by Corollary $\ref{Vcor}$ we have
\begin{equation*}\aligned\label{}\|\Phi u(t)-\Phi v(t)\|_{r,p}
&=\|G (P^\alpha V^\alpha u-P^\alpha V^\alpha(v))\|_{r,p}
\\ &\leq C\|V^\alpha u-V^\alpha (v)\|_{2a;k-b',\bar{c}}
\\ &\leq C(\|u\|_{a;k,c}+\|v\|_{a;k,c})\|u-v\|_{a;k,c}
\\ &\leq CM\|u-v\|_{a;k,c},
\endaligned
\end{equation*}
and similarly we have
\begin{equation*}\aligned\label{}\|\Phi u(t)-\Phi v(t)\|_{a;k,c}&=\|G (V^\alpha u-V^\alpha(v))\|_{a;k,c}
\\ &\leq C\|V^\alpha u-V^\alpha (v)\|_{2a;k-b',\bar{c}}
\\ &\leq C(\|u\|_{a;k,c}+\|v\|_{a;k,c})\|u-v\|_{a;k,c}
\\ &\leq CM\|u-v\|_{a;k,c}.
\endaligned
\end{equation*}

So for a sufficiently small choice of $M$, we can choose a $T$ such
that $\Phi$ sends $E_{T,M}$ into itself and is a contraction on
$E_{T,M}$. So by the contraction mapping principle, we have a unique
fixed point $u\in E_{T,M}$ provided our parameters satisfy all the
requisite inequalities.  Combining and simplifying these
inequalities, and allowing $s'=k-2-b+b'$ to define $s'$, we get the
following list of restrictions on the parameters:
\begin{equation}\label{conditions}\aligned& 1< p\leq c<\infty
\\&s':=k-2-b+b'
\\&k\geq 0,~~ b'\geq 1,~~s'c<n
\\&0<2a=k-n/c-b<1
\\&0\leq s'\leq k
\\&1<\frac{nc}{2n-s'c}\leq p
\\& 1\geq b'-b
\\& 1\leq b'+\frac{n}{c}-s'<2
\\& 2-2b'+s'\leq\frac{n}{p}\leq 2-b'+s'.
\endaligned
\end{equation}

This is not optimal, because of the presence of the ``extra"
parameter $b'$. However, this version does make it easy to ascertain
certain bounds on the original parameters. For example, the second
and seventh conditions require that $b\geq 0$, which provides a
lower bound of $n/p$ on the regularity of our initial data.

To eliminate the extra parameter $b'$, we remark that the
conditions force $1\leq b'<2$, and our optimal case ($b=0$) requires
$b'=1$. So setting $b'=1$, we let $k=1+b+s'$ define $s'$, and our
list of conditions becomes
\begin{equation}\label{conditionsshort}\aligned& 1< p\leq c<\infty
\\& b\geq 0
\\&s':=k-1-b
\\&k\geq 1,~~ s'c<n
\\&0<2a=k-n/c-b<1
\\&1<\frac{nc}{2n-s'c}\leq p
\\& 0\leq \frac{n}{c}-s'<1
\\& s'\leq\frac{n}{p}\leq 1+s'.
\endaligned
\end{equation}

\subsection{Proof of Theorem $\ref{special case 1}$}\label{Improvements for special cases}
The local existence result established in Section $\ref{Proof of thm1}$ required the initial data to have regularity no less than $n/p$.  In this section we prove Theorem $\ref{special case 1}$, which improves on that result by choosing specific values for the parameters.  Specifically, we begin with initial data $u_0\in H^{s_1,p}(\mathbb{R}^n)$ and we look for a solution $u$ to the LANS equation in
\begin{equation*}u\in BC([0,T):H^{s_1,p}(\mathbb{R}^n))\cap \dot{C}^T_{a;s_2,p}.
\end{equation*}
Our goal is to minimize the value of $s_1$, so we will choose $s_1<1$ and $s_2\geq 1$.  Then, applying the standard contraction method outlined in the previous section, we only need appropriate estimates for
\begin{equation}\label{improv1}\|e^{t\triangle}u_0\|_{a;s_2,p},
\end{equation}

\begin{equation}\label{improv2}\|\int_0^t e^{(t-s)\triangle}(\div(u\tensor u))ds\|_{0;s_1,p}+\|\int_0^t e^{(t-s)\triangle}(\div(1-\triangle)^{-1} (\nabla u\nabla ))ds\|_{0;s_1,p},
\end{equation}
and
\begin{equation}\label{improv3}\|\int_0^t e^{(t-s)\triangle}(\div(u\tensor u))ds\|_{a;s_2,p}+\|\int_0^t e^{(t-s)\triangle}(\div(1-\triangle)^{-1} (\nabla u\nabla ))ds\|_{a;s_2,p}.
\end{equation}

For $(\ref{improv1})$, we have
\begin{equation*}\|e^{t\triangle}u_0\|_{s_2,p}\leq t^{-(s_2-s_1)/2}\|u_0\|_{s_1,p},
\end{equation*}
and this forces $a=(s_2-s_1)/2$.  In order to estimate $(\ref{improv2})$ and $(\ref{improv3})$, we consider $p=2$ and $p>n$ separately, starting with $p=2$.  For $p=2$, the first term in $(\ref{improv2})$ satisfies
\begin{equation*}\aligned &\|\int_0^t e^{(t-s)\triangle}(\div(1-\triangle)^{-1} (\nabla u\nabla ))ds\|_{0;s_1,2}
\\ \leq &\int_0^t |t-s|^{-(n/4-n/2)/2}\|\nabla u\nabla u\|_{s_1-1,4}
\\ \leq &\int_0^t |t-s|^{-(n/2)/2}\|u\|_{1,2}\|u\|_{1,2}
\\ \leq &\|u\|_{\bar{a};1,2}\|u\|_{\bar{a};1,2}\int_0^t |t-s|^{-(n/4)}s^{-(1-s_1)/2-(1-s_1)/2}ds
\\ \leq &\|u\|_{\bar{a};1,2}\|u\|_{\bar{a};1,2} t^{-n/4+s_1},
\endaligned
\end{equation*}
where we set $\bar{a}=(1-s_1)/2$.  For the estimates to apply, and for the exponent of $t$ to have the desired sign, the parameters are required to satisfy
\begin{equation*}\aligned n/2&<2,
\\ s_1&\geq n/4.
\endaligned
\end{equation*}
We observe that the first condition requires $n=3$, which means $s_1\geq 3/4$.  We also remark that a straightforward modification of the estimate will allow $n\geq 4$, provided $s_1$ is no longer required to be less than one.  Estimating the second term in $(\ref{improv2})$, we have
\begin{equation*}\aligned &\int_0^t |t-s|^{-(3/4+n/\tilde{p}-n/2)/2}\|\div(u\tensor u)\|_{0,\tilde{p}}ds
\leq \int_0^t |t-s|^{-(3/4+n/\tilde{p}-n/2)/2}\|u\|_{1,2}\|u\|_{\bar{p}} ds
\\ &\leq \int_0^t |t-s|^{-(3/4+n/\bar{p})/2}\|u\|_{1,2}\|u\|_{1,2} ds \leq \|u\|_{\bar{a};1,2}\int_0^t |t-s|^{-(3/4+n/2-1)/2}s^{-(1-3/4)}ds
\\ &\leq \|u\|_{\bar{a};1,2}t^{-(1-3/4+n/2)/2+1}=\|u\|_{\bar{a};s_2,2}t^{-n/4+7/8}.
\endaligned
\end{equation*}

Observing that for $n=3$, $-n/4+7/8$ has the desired sign, the requirements of the parameters are
\begin{equation*}\aligned 1/\tilde{p}&=1/2+1/\bar{p}
\\ \bar{p}&=\frac{2n}{n-2s_2}
\\ (3/4+n/2-1)/2&=(n/2-1/4)/2<1.
\endaligned
\end{equation*}

For $(\ref{improv3})$, we have
\begin{equation*}\aligned &\|\int_0^t e^{(t-s)\triangle}(\div(1-\triangle)^{-1} (\nabla u\nabla u))ds\|_{a;s_2,2}
\\ \leq &t^a\int_0^t |t-s|^{-(n/4)}\|\nabla u\nabla u\|_{1-1,4} \leq t^a\int_0^t |t-s|^{-(n/4)}\|u\|^2_{1,2}
\\ \leq &t^a\|u\|^2_{a;1,2}\int_0^t |t-s|^{-(n/4)}s^{-(1-s_1)}ds \leq \|u\|^2 t^{-(n/4+1-s_1)+1+(1-3/4)/2+a},
\endaligned
\end{equation*}

and
\begin{equation*}\aligned &\|\int_0^t e^{(t-s)\triangle}(\div(u\tensor u)ds\|_{a;s_2,2} \leq t^a\int_0^t |t-s|^{-(1+n/\tilde{p}-n/2)/2}\|\div(u\tensor u)\|_{0,\tilde{p}}ds
\\ \leq &t^a\int_0^t |t-s|^{-(1+n/\tilde{p}-n/2)/2}\|u\|_{1,2}\|u\|_{\bar{p}} ds \leq t^a\int_0^t |t-s|^{-(1+n/\bar{p})/2}\|u\|_{1,2}\|u\|_{1,2} ds
\\ \leq &\|u\|_{a;1,2}t^a\int_0^t |t-s|^{-(1+n/2-1)/2}s^{-(1-3/4)}ds \leq \|u\|_{a;1,2}t^{-(n/2)/2+1+(1-3/4)/2+a}.
\endaligned
\end{equation*}

So choosing $s_1=3/4$, $n=3$, and $s_2=1$ satisfies all the requirements on the parameters, which finishes the proof for the case $p=2$.

For $p>n$, the first term of $(\ref{improv2})$ satisfies
\begin{equation*}\aligned &\|\int_0^t e^{(t-s)\triangle}(\div(1-\triangle)^{-1} (\nabla u\nabla ))ds\|_{0;s_1,p}
 \leq \int_0^t |t-s|^{-(n/2p)}\|\nabla u\nabla u\|_{s_1-1,2p}
\\ \leq &\int_0^t |t-s|^{-(n/2p)}\|u\|^2_{1,p} \leq \|u\|^2_{a;1,p}\int_0^t |t-s|^{-(n/2p)}s^{-(1-s_1)}ds
\\ \leq &\|u\|^2 t^{-n/2p+s_1}.
\endaligned
\end{equation*}

Here, the principle restriction is that $s_1\geq n/2p$.  For the second term, we have
\begin{equation*}\aligned \int_0^t |t-s|^{-(s_1+n/p)/2}\|\div(u\tensor u)\|_{0,2p}ds
 &\leq \int_0^t |t-s|^{-(s_1+n/p)/2}\|u\|_{1,p}\|u\|_{0,p} ds
\\ &\leq \|u\|_{a;1,p}\|u\|_{0;0,p}\int_0^t |t-s|^{-(s_1+n/p)/2}s^{-(1-s_1)/2} ds
\\ &\leq \|u\|_{a;1,p}\|u\|_{0;0,p}t^{-(1+s_1)/2+1}.
\endaligned
\end{equation*}

Since $s_1<1$, the exponent on $t$ in the last line has the desired sign.  For $(\ref{improv3})$, we have
\begin{equation*}\aligned &\|\int_0^t e^{(t-s)\triangle}(\div(1-\triangle)^{-1} (\nabla u\nabla ))ds\|_{a;k,p}
\leq t^a\int_0^t |t-s|^{-(n/2p)}\|\nabla u\nabla u\|_{1-1,2p}
\\ \leq &t^a\int_0^t |t-s|^{-(n/2p)}\|u\|^2_{1,p} \leq t^a\|u\|^2_{a;1,p}\int_0^t |t-s|^{-(n/2p)}s^{-(1-s_1)}ds
\\ \leq &\|u\|^2 t^{-n/2p+1-(1-s_1)/2},
\endaligned
\end{equation*}

and
\begin{equation*}\aligned &t^a\int_0^t |t-s|^{-(1+n/p)/2}\|\div(u\tensor u)\|_{0,2p}ds
 \leq t^a\int_0^t |t-s|^{-(1+n/p)/2}\|u\|_{1,p}\|u\|_{0,p} ds
\\ &\leq \|u\|_{a;1,p}\|u\|_{0;0,p}t^a \int_0^t |t-s|^{-(1+n/p)/2}s^{-(1-s_1)/2} \leq \|u\|_{a;1,p}^2t^{-(1+n/p)/2+1}.
\endaligned
\end{equation*}

This finishes the proof of Theorem $\ref{special case 1}$.

\section{Solutions in $L^a((0,T):H^{s,p}(\mathbb{R}^n))$}\label{Solutions in the class of integral norms in time}

In this section we address Theorems $\ref{special case 2}$ and $\ref{thm2}$. We begin in Section $\ref{Preliminary results}$ with necessary supporting results.  Section $\ref{Proof of thm2}$ contains the proofs of the two theorems.

\subsection{Preliminary results}\label{Preliminary results}
Our first supporting result is Lemma $3.2$ in \cite{KP} and involves the operator
$\Gamma$.
\begin{proposition}\label{gamma prop 2}Let $1<q_0\leq q_1<\infty$, $s_0\leq s_1$, and assume
$0<(s_1-s_0+n/q_0-n/q_1)/2=1/\sigma \leq 1/q_0$.  Then $\Gamma$ maps
$H^{s_0,q_0}(\mathbb{R}^n)$ continuously into $L^\sigma((0,\infty):H^{s_1,q_1}(\mathbb{R}^n))$,
with the estimate
\begin{equation*}\(\int_0^\infty \|\Gamma u\|_{H^{s_1,q_1}(\mathbb{R}^n)}^\sigma\)^{1/\sigma}
\leq C\|u\|_{H^{s_0,q_0}(\mathbb{R}^n)}.
\end{equation*}
\end{proposition}

\begin{proof}We first observe that
$(s_1-s_0+n/q_0-n/q_1)/2=1/\sigma< 1$ implies $s_1-s_0< 2$.  So
without loss of generality, we assume $s_0=0$ and $s_1\in [0,2)$.  Next, we define the quasi-linear operator $K$ by
\begin{equation*}(Kf)(t)=\|e^{t\triangle}f\|_{H^{s_1,q_1}(\mathbb{R}^n)},
\end{equation*}
where $s_1$ and $q_1$ have been fixed.  Using the heat kernel estimate, we have that
\begin{equation}\label{chi}(Kf)(t)=\|e^{t\triangle}f\|_{s_1,q_1}\leq C
t^{-1/\sigma}\|f\|_{L^{q_0}}=C(Tf)(t),
\end{equation}
so $K$ is a quasi-linear map that
satisfies
\begin{equation*}K:L^{q_0}(\mathbb{R}^n)\rightarrow L_{1/\sigma,\infty}(I).
\end{equation*}
The Proposition then follows from an application of the Marcinkiewitz interpolation theorem.
\end{proof}

The following corollary follows from an application of the dominated convergence theorem.
\begin{corollary}\label{gamma prop 2 cor}For any $\varepsilon>0$, there exists a
$T$ which depends only on $\varepsilon$ and $\|u\|_{s_0,q_0}$
such that
\begin{equation*}\(\int_0^T \|e^{t\triangle}u\|_{s_1,q_1}^\sigma
dt\)^{1/\sigma}\leq \varepsilon,
\end{equation*}
for all $0<t<T$.
\end{corollary}

Next, we consider the operator $V^\alpha$ on our integral norm
space.
\begin{proposition}\label{Vprop2}Let $b\geq 1$, $1<q,p<\infty$, with $2/p-1/q<1$ and $0\leq n(2q-p)/pq\leq b-1$.  Then
\begin{equation*}V^\alpha:L^{\sigma}((0,T):H^{b,p}(\mathbb{R}^n))\rightarrow L^{\sigma/2}((0,T):H^{b-1,q}(\mathbb{R}^n))
\end{equation*}
with the estimate
\begin{equation}\label{Vest2}\(\int_0^T \|V^\alpha(u(s))\|^{\sigma/2}_{b-1,q}ds\)^{2/\sigma}\leq
\(\int_0^T \|u(s)\|^{\sigma}_{b,p}ds\)^{2/\sigma}.
\end{equation}
\end{proposition}
This follows directly from Proposition $\ref{Vprop}$.  We also have
that
\begin{equation*}\aligned
&\(\int_0^T\|V^\alpha(u(s))-V^\alpha(v(s))\|^{\sigma/2}_{b-1,q}ds\)^{2/\sigma}
\\\leq&\(\int_0^T(\|v(s)\|_{b,p}+\|u(s)\|_{b,p})^{\sigma/2}
(\|v(s)-u(s)\|_{b,p})^{\sigma/2}ds\)^{2/\sigma}
\\\leq&\(\int_0^T(\|v(s)\|_{b,p}+\|u(s)\|_{b,p})^{\sigma}ds\)^{2/\sigma}
\(\int_0^T\|v(s)-u(s)\|_{b,p}^{\sigma}ds\)^{2/\sigma}
\endaligned
\end{equation*}
where we used Holder's inequality and Minkowski's inequality.  This
gives an analog to Corollary $\ref{Vcor}$.
\begin{corollary}\label{Vcor2}With the same assumptions on the parameters as in Proposition
$\ref{Vprop2}$, we have that if $u,v\in L^{\sigma}((0,T):H^{b,p}(\mathbb{R}^n))$
then

\begin{equation}\label{Vcont2}\aligned
&\(\int_0^T\|V^\alpha(u(s))-V^\alpha(v(s))\|^{\sigma/2}_{b-1,q}ds\)^{2/\sigma}
\\\leq&\(\int_0^T(\|v(s)\|_{b,p}+\|u(s)\|_{b,p})^{\sigma}ds\)^{2/\sigma}
\(\int_0^T\|v(s)-u(s)\|_{b,p}^{\sigma}ds\)^{2/\sigma}.
\endaligned
\end{equation}
\end{corollary}

Our next set of results involve the operator $G$.
\begin{proposition}\label{G prop 1}Let $1\leq q'\leq q''<\infty$, $s'\leq s''$,
$1<\sigma'<\sigma''<\infty$, and let
$1/\sigma'-1/\sigma''=1-(s''-s'+n/q'-n/q'')/2.$  Then for any $T\in
(0,\infty]$, $G$ maps $L^{\sigma'}((0,T):H^{s',q'}(\mathbb{R}^n))$ continuously
into $L^{\sigma''}((0,T):H^{s'',q''}(\mathbb{R}^n)).$
\end{proposition}

Using Proposition $\ref{G}$, we observe that
\begin{equation*}\|G u(t)\|_{s'',q''}\leq C \int_0^T |t-s|^{(s''-s'+n/q'-n/q'')/2} \|u(s)\|_{s',q'}ds.
\end{equation*}
Applying the Hardy-Littlewood-Sobolev Theorem, with $1/r=(s''-s'+n/p'-n/p'')/2$, we have
\begin{equation*}\aligned (\int_0^T \|G u\|_{s'',q''}^{\sigma''}dt)^{1/\sigma''}
\leq C\|f\|_{L^{\sigma'}(I)}
\endaligned
\end{equation*}
where $f(t)=\|u(t)\|_{s',q'}$.  This completes the proof.

Our next result also involves the operator $G$.

\begin{proposition}\label{G prop 2}Let $1<q'\leq q''<\infty$, $s'\leq s''$ and assume
$1/q''\leq 1/\sigma=1-(s''-s'+n/q'-n/q'')/2\leq 1$.  Then $G$ maps
$L^{\sigma}((0,T):H^{s',q'}(\mathbb{R}^n))$ continuously into
$BC([0,T):H^{s'',q''}(\mathbb{R}^n))$.
\end{proposition}

The proof is a duality argument, which we include for completeness.  We begin by defining $H$ by $(Hf)(s,x)=e^{s\triangle}f(s,x)$.  We recall that the dual of the space $L^{q}((0,T):L^{p})$ is
$L^{\bar{q}}((0,T):L^{\bar{p}})$ (where $\bar{q}$ and $\bar{p}$
denote the conjugate exponents to $q$ and $p$ and $1\leq
p,q<\infty$).  So for any $g\in
L^{\infty}((0,T):H^{-s'',\bar{q}''}(\mathbb{R}^n))$, we have
\begin{equation}\label{stepf}\aligned&\int_0^T \int_{\mathbb{R}^n}
g(s,x)e^{s\triangle}(1-\triangle)^{s''/2}f(s,x)dxds
\\\leq &C\int_0^T\int_{\mathbb{R}^n}e^{s\triangle}(1-\triangle)^{(s''-s')/2}
g(s,x)(1-\triangle)^{s'/2}f(s,x)dxds
\\\leq &C\int_0^T\|e^{s\triangle}g(s)\|_{H^{s''-s',\bar{q}'}(\mathbb{R}^n)} \|f(s)\|_{s',q'}ds
\\\leq &C(\int_0^T
\|e^{s\triangle}g(s)\|^{\bar{\sigma}}_{s''-s',\bar{q}'}ds)^{1/\bar{\sigma}}
(\int_0^T \|f(s)\|^{\sigma}_{s',q'}ds)^{1/\sigma}
\\\leq
&C\sup_s\|g(s)\|_{\bar{q}''}(\int_0^T\|f(s)\|^{\sigma}_{s',q'}ds)^{1/\sigma},
\endaligned
\end{equation}
where the last line is a slight generalization of Proposition
$\ref{gamma prop 2}$.  Since $g$ is an arbitrary element of the dual
space of $L^1((0,T):H^{s'',q''}(\mathbb{R}^n))$, we have
\begin{equation}\label{G step 2}\int_0^T \|e^{s\triangle}f(s)\|_{s'',q''}ds
\leq C(\int_0^T \|f(s)\|^{\sigma}_{s',q'}ds)^{1/\sigma}.
\end{equation}

To finish the proposition, making liberal use of the change of
variables formula and using $(\ref{G step 2})$, we have
\begin{equation}\label{last step}\aligned \|G\cdot f(t)\|_{s'',q''}&=\|\int_0^te^{(t-s)\triangle}f(s)ds\|_{s'',q''}
\\&\leq \int_0^t\|e^{(t-s)\triangle}f(s)\|_{s'',q''}ds
\\&=\int_0^t\|e^{s\triangle}f(t-s)\|_{s'',q''}ds
\\ &\leq C(\int_0^t\|f(t-s)\|^{\sigma}_{s',q'}ds)^{1/\sigma}
\\ &\leq C(\int_0^T\|f(s)\|^{\sigma}_{s',q'}ds)^{1/\sigma},
\endaligned
\end{equation}
which proves the proposition.
\subsection{Proof of Theorems $\ref{special case 2}$ and $\ref{thm2}$}\label{Proof of thm2} To prove Theorem $\ref{thm2}$, as in
Section $\ref{Proof of thm1}$, we begin with the nonlinear map
\begin{equation*}\Phi u=\Gamma u_0-G\cdot P^\alpha(V^\alpha(u))
\end{equation*}
and the space $F_{T,M}$ defined to be the space of all
\begin{equation*}v\in BC([0,T):H^{r,p}) \cap L^a((0,T):H^{k,c})
\end{equation*}
such that
\begin{equation*}\sup_{0\leq t\leq T}\|v(t)-\Gamma u_0\|_{r,p}
+\(\int_0^T \|v(s)\|^{a}_{k,c}ds\)^{1/a}\leq M.
\end{equation*}

Using the same argument used in Section $\ref{Proof of thm1}$, we
get that $\Phi$ will be a contraction mapping provided the following
list of conditions is satisfied:
\begin{equation}\label{conditions2}\aligned& 1< p\leq c<\infty
\\&s':=k-2+b'-b
\\&k\geq 1,~~ b'\geq 1,~~s'c<n
\\&0<2/a=k-n/c-b<1
\\&0\leq s'\leq k-1
\\&1<=\frac{nc}{2n-s'c}
\\& 1\geq b'-b
\\& k-b'\leq \frac{n}{p}+b\leq k
\\& a/2\leq p\leq a.
\endaligned
\end{equation}

We observe that as in the previous case, these conditions require
that $b\geq 0$.  We also record the simplified list that arises from
setting $b'=1$:

\begin{equation}\label{conditions2short}\aligned& 1< p\leq c<\infty
\\ &b\geq 0
\\&s':=k-1-b
\\&k\geq 1,~~ s'c<n
\\&0<2/a=k-n/c-b<1
\\& k-1\leq \frac{n}{p}+b\leq k
\\& a/2\leq p\leq a.
\endaligned
\end{equation}

In order to prove $\ref{special case 2}$, we employ an argument similar to the one given for the proof of Theorem $\ref{special case 1}$ in Section $\ref{Improvements for special cases}$.

\section{Uniform-in-time bounds for solutions to the LANS equation}\label{Global
Existence in Sobolev space}
In this section we prove Theorem $\ref{globathm2}$.  The proof has two main steps.  First, we prove that  $\|u(t,\cdot)\|_{1,2}$ is a decreasing function of time.  Then we bound $\|u(t,\cdot)\|_{2,2}$ by a function of $\|u(t,\cdot)\|_{1,2}$, and thus by the first step have our desired bound on $\|u(t,\cdot)\|_{2,2}$.  We remark that higher regularity bounds can be obtained by successive steps.  For additional details, see \cite{MS}.

We also note that, when applying this result to our local existence theorems, only the first step is required to control the solutions from Theorems $\ref{special case 1}$ and $\ref{special case 2}$.  Finally, for notational simplicity, we set $\|u(t,\cdot)\|_{r,p}=\|u(t)\|_{r,p}$.

\begin{proof}
We start by recalling $(\ref{LANS})$
\begin{equation*}\partial_t u+(u\cdot\nabla)u+\div
\tau^\alpha u=-(1-\alpha^2\triangle)^{-1}\nabla p+\nu\triangle
\end{equation*}
and stating an equivalent form (see Section $3$ of \cite{MS})
\begin{equation}\label{LANSglobalpi}\aligned &\partial_t (1-\alpha^2\triangle)u
+\nabla_u[(1-\alpha^2\triangle)u] -\alpha^2(\nabla
u)^T\cdot\triangle u
\\=&-(1-\alpha^2\triangle)A u-\nabla p.
\endaligned
\end{equation}

%start of H^1 bound
To start, we take the $L^2$ product of $(\ref{LANSglobalpi})$ with
$u$. We get
\begin{equation*}I_1+I_2+I_3=J_1+J_2
\end{equation*}
where
\begin{equation*}\aligned I_1&=(\partial_t (1-\alpha^2\triangle)u,u)
\\ I_2&=(\nabla_u u,u)
\\ I_3&=-\alpha^2\((\nabla_u \triangle u,u)+((\nabla u)^T\cdot \triangle u, u)\)
\\ J_1&=-((1-\alpha^2 \triangle)(A u),u)
\\ J_2&=(\nabla p,u).
\endaligned
\end{equation*}

We start with $I_1$, which becomes
\begin{equation*}\aligned I_1&=(\partial_t u,u)-\alpha^2(\triangle \partial_t u,u)
\\ &=\frac{1}{2}\partial_t (\|u\|_{L^2}^2+\alpha^2\|A^{1/2}u\|^2_{L^2}).
\endaligned
\end{equation*}
Applying integration by parts to $I_2$, $I_3$, and $J_2$ and recalling that $\div u=0$, we get that all three terms vanish.  For $J_1$ we have 
\begin{equation*}J_1=-((1-\alpha^2 \triangle)(A u),u)=-(A^{1/2}u,A^{1/2}u)-\alpha^2(Au,Au).
\end{equation*}

Applying this to $(\ref{LANSglobalpi})$, we get
\begin{equation}\label{global11}\frac{1}{2}\partial_t (\|u(t)\|^2_{L^2}+\alpha^2\|u(t)\|^2_{\dot{H}^{1,2}})
\leq -(\|A^{1/2}u(t)\|^2_{L^2}+\alpha^2\|Au(t)\|_{L^2}^2),
\end{equation}
where $\dot{H}$ denotes the homogeneous Sobolev norm.  This proves
that $\|u(t)\|_{H^{1,2}}$ is decreasing in time.

%end H^1 bound, start of H^2 est.
We remark that this is sufficient to extend the local solutions of Theorems $\ref{special case 1}$ and $\ref{special case 2}$.  However, in order to extend the local solutions of Theorems $\ref{thm1}$ and $\ref{thm2}$, we require the following additional estimate.

Applying $A$ to $(\ref{LANS})$ and taking the $L^2$ product with $Au$, we get
\begin{equation}\label{aldo}(\partial_t Au,Au)+(A^2u,Au)+(AP^\alpha(\nabla_u u+\div\tau^\alpha u,Au))=0.
\end{equation}
The first piece satisfies
\begin{equation}\label{aldo1}(\partial_t Au,Au)=\frac{1}{2}\partial_t \|Au\|^2_{L^2}
\end{equation}
and the second satisfies
\begin{equation}\label{aldo2}(A^2 u,Au)=(A^{3/2}u,A^{3/2})=\|A^{3/2}\|^2_{L^2}.
\end{equation}

%supporting inequalities
To handle the last term of $(\ref{aldo})$, we write it as
\begin{equation*}(AP^\alpha(\nabla_u u),Au)+(AP^\alpha\div\tau^\alpha u,Au)=K_1+K_2.
\end{equation*}

To proceed, we will need two inequalities. The first is the well
known Sobolev embedding:
\begin{equation*}\|u\|_{L^\infty}\leq C\|u\|_{H^{k,2}}
\end{equation*}
provided $2k>3$.  The second is called a Ladyzhenskaya inequality,
and is $(5.3)$ in \cite{MS}:
\begin{equation}\label{AgLa}\aligned
\|u\|_{\dot{H}^i}&\leq C\|u\|^{1-i/m}_{L^2}\|u\|^{i/m}_{\dot{H}^m},
\endaligned
\end{equation}
where $\dot{H}$ is the homogeneous Sobolev space.
%back to H^2 est.

Starting with $K_1$, we have
\begin{equation}\label{glob1d}\aligned (AP^\alpha(\nabla_u u),Au) &=(A^{1/2}(\nabla u\cdot u),A^{3/2}u)
\\ &\leq C\|A^{3/2}u\|_{L^2}(\|(A^{1/2}\nabla u)u\|_{L^2}+\|(A^{1/2}u)\nabla u\|_{L^2})
\\ &\leq C\|u\|_{\dot{H}^3} (\|u\|_{L^\infty} \|A^{1/2}\nabla u\|_{L^2}+\|A^{1/2}u\|_{L^\infty}\|\nabla u\|_{L^2})
\\ &\leq C\|u\|_{\dot{H}^3}
\(\|u\|_{L^\infty}\|u\|_{\dot{H}^2}+\|A^{1/2}u\|_{L^\infty}\|u\|_{\dot{H}^1}\).
\endaligned
\end{equation}
By Sobolev embedding and the observation that $\|u\|_{H^{s,p}}\leq \|u\|_{L^p}+\|u\|_{\dot{H}^{s,p}}$, 
we have
\begin{equation}\label{agla2}\aligned \|u\|_{L^\infty}&\leq
C\|u\|_{H^{k_1}}\leq C(\|u\|_{L^2}+\|u\|_{\dot{H}^{k_1}}) \leq
C(\|u\|_{H^1}+\|u\|_{\dot{H}^{k_1}})
\\ \|A^{1/2}u\|_{L^\infty}&\leq C\|u\|_{H^{k_2}}\leq
C(\|\nabla u\|_{L^2}+\|u\|_{\dot{H}^{k_2}})\leq
C(\|u\|_{H^1}+\|u\|_{\dot{H}^{k_2}})
\endaligned
\end{equation}
where $k_1=3/2+\varepsilon$ and $k_2=5/2+\delta$ for positive
numbers $\varepsilon$ and $\delta$.  So $(\ref{glob1d})$ becomes
\begin{equation}\label{glob1c}\aligned (AP^\alpha(\nabla_u u),Au)
&\leq C\|u\|_{\dot{H}^3}\|u\|_{\dot{H}^2}\|u\|_{H^1}
+C\|u\|_{\dot{H}^3} \|u\|_{\dot{H}^2}\|u\|_{\dot{H}^{k_1}}
\\&+ C\|u\|_{\dot{H}^3} \|u\|_{H^1}\|u\|_{\dot{H}^{k_2}}+C\|u\|_{\dot{H}^3}
\|u\|_{H^1}^2.
\endaligned
\end{equation}

By $(\ref{AgLa})$, we have
\begin{equation}\aligned \label{agla3}\|u\|_{\dot{H}^2}&=\|\nabla u\|_{\dot{H}^1}
\leq C\|\nabla u\|_{L^2}^{1/2}\|\nabla u\|_{\dot{H}^2}^{1/2} \leq
C\|u\|_{\dot{H}^1}^{1/2}\|u\|_{\dot{H}^3}^{1/2}
\\ \|u\|_{\dot{H}^{k_1}}&=\|\nabla u\|_{\dot{H}^{k_1-1}} \leq
C\|u\|^{1-(k_1-1)/2}_{\dot{H}^{1}} \|u\|^{(k_1-1)/2}_{\dot{H}^{3}}
\\ \|u\|_{\dot{H}^{k_2}}&=\|\nabla u\|_{\dot{H}^{k_2-1}} \leq
C\|u\|^{1-(k_2-1)/2}_{\dot{H}^{1}} \|u\|^{(k_2-1)/2}_{\dot{H}^{3}}.
\endaligned
\end{equation}

Applying $(\ref{agla3})$ to $(\ref{glob1c})$, we have
\begin{equation}\label{glob1b}\aligned (AP^\alpha(\nabla_u u),Au)
\leq &C\|u\|_{\dot{H}^3}^{1+k_1/2}\|u\|_{\dot{H}^1}^{2-k_1/2}+
C\|u\|_{\dot{H}^3}^{(k_2+1)/2}\|u\|_{\dot{H}^1}^{(k_2+3)/2}
\\ &+C\|u\|_{\dot{H}^3}^{3/2}\|u\|_{H^1}^{3/2}+C\|u\|_{\dot{H}^3}\|u\|_{H^1}^2.
\endaligned
\end{equation}

Choosing $\varepsilon=\delta=1/4$, we get
\begin{equation}\aligned \label{glob1a}(AP^\alpha(\nabla_u u),Au)\leq
&C\|u\|^{15/8}_{\dot{H}^3}(\|u\|_{\dot{H}^1}^{9/8}+\|u\|_{\dot{H}^1}^{23/8})
\\ +&C\|u\|_{\dot{H}^3}^{3/2}\|u\|_{H^1}^{3/2}
+C\|u\|_{\dot{H}^3}\|u\|_{H^1}^2,
\endaligned
\end{equation}
which finishes our $K_1$ estimate.  For $K_2$, we have
\begin{equation}\label{glob1er}(A(\div\tau^\alpha)(u),Au)
\leq \|u\|_{\dot{H}^2}\|A(\div\tau^\alpha)(u)\|_{L^2}.
\end{equation}

To estimate the second term, we remark that it is sufficient to
consider $A(1-\alpha^2\triangle)^{-1}\div(\nabla u\cdot\nabla u)$,
and we have
\begin{equation*}\aligned \|A(1-\alpha^2\triangle)^{-1}\div(\nabla u\cdot\nabla u)\|_{L^2}
&\leq \|\div(\nabla u\cdot\nabla u)\|_{L^2}
\\ &\leq C\|\nabla u\|_{L^\infty}\|u\|_{\dot{H}^2}.
\endaligned
\end{equation*}
Plugging this back into $(\ref{glob1er})$ and using $(\ref{agla2})$
and $(\ref{agla3})$ gives
\begin{equation}\label{glob1b2}\aligned (A(\div\tau^\alpha)(u),Au)
&\leq C\|u\|_{\dot{H}^2}^2 \|\nabla u\|_{L^\infty}
\\ &\leq C(\|u\|_{\dot{H}^2}^2\|u\|_{H^1}
+ \|u\|_{\dot{H}^2}^2 \|u\|_{\dot{H}^{k_2}})
\\ &\leq C(\|u\|_{\dot{H}^3} \|u\|_{H^1}^2
+ \|u\|_{\dot{H}^3}^{15/8}\|u\|_{H^1}^{23/8}).
\endaligned
\end{equation}

Combining $(\ref{glob1a})$ and $(\ref{glob1b2})$ gives
\begin{equation}\label{aldo16}\aligned (AP^\alpha(\nabla_u u+\div\tau^\alpha u),Au)
\leq &C\|u\|^{15/8}_{\dot{H}^3}(\|u\|_{\dot{H}^1}^{9/8}
+\|u\|_{\dot{H}^1}^{23/8})
\\ +&C\|u\|_{\dot{H}^3}^{3/2}\|u\|_{H^1}^{3/2}
+C\|u\|_{\dot{H}^3}\|u\|_{H^1}^2.
\endaligned
\end{equation}
Applying Young's multiplicative inequality with $p=16/15$ and $p'=16$, we get
\begin{equation}\label{aldo4}\|u\|^{15/8}_{\dot{H}^3}(\|u\|_{\dot{H}^1}^{9/8}
+\|u\|_{\dot{H}^1}^{23/8}) \leq C\varepsilon\|u\|_{\dot{H}^3}^2
+\frac{C}{\varepsilon}(\|u\|_{H^1}^{18}+\|u\|_{H^1}^{46}).
\end{equation}

Choosing $\varepsilon=(4C)^{-1}$, $(\ref{aldo4})$ becomes
\begin{equation}\label{aldo17}\|u\|^{15/8}_{\dot{H}^3}(\|u\|_{\dot{H}^1}^{9/8}
+\|u\|_{\dot{H}^1}^{23/8}) \leq \frac{1}{4}\|u\|_{\dot{H}^3}^2
+C(\|u\|_{H^1}^{18}+\|u\|_{H^1}^{46}).
\end{equation}

Similarly, Young's multiplicative inequality applied to the second term on the right hand side of $(\ref{aldo16})$ gives
\begin{equation}\label{aldo18}\|u\|_{\dot{H}^3}^{3/2}\|u\|_{H^1}^{3/2}
+C\|u\|_{\dot{H}^3}\|u\|_{H^1}^2\leq \frac{1}{4}\|u\|_{\dot{H}^3}^2
+ C(\|u\|_{H^1}^6+\|u\|_{H^1}^4).
\end{equation}

Using $(\ref{aldo17})$ and $(\ref{aldo18})$ in $(\ref{aldo16})$
gives
\begin{equation}\label{aldo3}(AP^\alpha(\nabla_u u+\div\tau^\alpha u),Au)
\leq \frac{1}{2}\|u\|_{\dot{H}^3}^2
+C(\|u\|_{H^1}^{18}+\|u\|_{H^1}^{46} +\|u\|_{H^1}^6+\|u\|_{H^1}^4).
\end{equation}

Finally, using $(\ref{aldo1})$, $(\ref{aldo2})$ and $(\ref{aldo3})$
in $(\ref{aldo})$ gives
\begin{equation*}\aligned \frac{1}{2}\partial_t \|u(t)\|^2_{\dot{H}^2}
&\leq \frac{-1}{2}\|u(t)\|^2_{\dot{H}^3}
+C(\|u(t)\|_{H^1}^{18}+\|u(t)\|_{H^1}^{46}+\|u\|_{H^1}^6+\|u\|_{H^1}^4)
\\ &\leq \frac{-1}{2}\|u(t)\|_{\dot{H}^3}^2
+C(\|u_0\|_{H^1}^{18}+\|u_0\|_{H^1}^{46}
+\|u_0\|_{H^1}^6+\|u_0\|_{H^1}^4),
\endaligned
\end{equation*}
where the last line used $(\ref{global11})$.  So, for any $t$ such
that
\begin{equation*}\|u(t)\|_{\dot{H}^3}\geq C(\|u_0\|_{H^1}^{18}
+\|u_0\|_{H^1}^{46}+\|u_0\|_{H^1}^6+\|u_0\|_{H^1}^4)^{1/2},
\end{equation*}
we get that $\|u(t)\|_{\dot{H}^2}$ is decreasing as a function of
time at $t$.  So our last task is to show that $\|u\|_{\dot{H}^2}$
is bounded provided
\begin{equation*}\|u(t)\|_{\dot{H}^3}< C(\|u_0\|_{H^1}^{18}
+\|u_0\|_{H^1}^{46}+\|u_0\|_{H^1}^6+\|u_0\|_{H^1}^4)^{1/2}.
\end{equation*}
%end of bulk of H^2, still need silly trick

To handle this case, we again use $(\ref{AgLa})$, and get
\begin{equation}\label{cmpunk}\aligned \|u(t)\|_{\dot{H}^2}&\leq C\|u(t)\|_{L^2}^{1/3}\|u(t)\|_{\dot{H}^3}^{2/3}
\\ &\leq C\|u_0\|_{L^2}^{1/3}(\|u_0\|_{H^1}^{10}
+\|u_0\|_{H^1}^{2}+\|u_0\|_{H^1}^6+\|u_0\|_{H^1}^4)^{1/3}.
\endaligned
\end{equation}
Combining $(\ref{global11})$ and $(\ref{cmpunk})$, we get that 
\begin{equation*}\|u(t)\|_{2,2}\leq CM,
\end{equation*}
where $M$ depends only on $\|u_0\|_{1,p}$, which proves Theorem $\ref{globathm2}$.
\end{proof}

\section{Higher regularity for the local existence result}\label{Higher regularity for the local existence result}
In this section we exploit the regularizing effect of the heat kernel to prove that the local solutions constructed in Theorems $\ref{special case 1}$ through $\ref{thm2}$ possess more regularity than mandated by the theorem.  Our results also quantify the blow-up that occurs in these higher regularity norms as $t\rightarrow 0$.  We use an induction argument inspired by the results in \cite{katoinduction} for the Navier-Stokes equation.    

This additional regularity also allows us to apply the $\emph{a~priori}$ estimate from Theorem $\ref{globathm2}$ to our local existence results.  

\begin{proposition}\label{higher regularity theorem}Let $u_0\in H^{s_1,p}(\mathbb{R}^n)$ be divergence-free, and let $s_1<1$.  Let $u$ be a solution to the LANS equation ($\ref{LANS})$ such that   
\begin{equation*}u\in BC([0,T):H^{s_1,p}(\mathbb{R}^n))\cap \dot{C}^T_{(s_2-s_1)/2;s_2,p},
\end{equation*}
where $0<s_2-s_1<1$ and $s_2\geq 1$.  Then for all $r\geq s_2$, we have that $u\in \dot{C}^T_{(r-s_1)/2;r,p}$.
\end{proposition}

We have an analogous result for the $L^a((0,T):H^{s,p}(\mathbb{R}^n))$ case.

\begin{proposition}\label{higher regularity theorem2}Let $k>s_2>s_1$, with $s_2\geq 1$, and let $\varepsilon$ be a small positive number.  Then, for $k-s_2=s_2-s_1=\varepsilon$, for any solution $u$ to the LANS equation $(\ref{LANS})$ where 
\begin{equation*}u\in BC([0,T):H^{s_1,p}(\mathbb{R}^n)\cap L^{2/(s_2-s_1)}((0,T):H^{s_2,p}(\mathbb{R}^n)),
\end{equation*}
we have that $u\in L^1((0,T):H^{k,p}(\mathbb{R}^n))$ provided $s_2\geq n/p$.
\end{proposition}

Before addressing the proofs, we remark that both propositions can be readily extended to additional cases.  In Proposition $\ref{higher regularity theorem}$, the requirement that $s_1<1$ can be removed and the assumption that $u\in BC([0,T):H^{s_1,p}(\mathbb{R}^n))\cap \dot{C}^T_{(s_2-s_1)/2;s_2,p}$ can be changed to $u\in BC([0,T):H^{s_1,p}(\mathbb{R}^n))\cap \dot{C}^T_{a;s_2,q}$ for $q>p$.  There are similar extensions for Proposition $\ref{higher regularity theorem2}$, with the additional requirement that, for $k$ such that $k-s_1\geq 2$, $u\in L^a((0,T):H^{k,q}(\mathbb{R}^n)$ requires $q>p$.   

The proofs of Propositions $\ref{higher regularity theorem}$ and $\ref{higher regularity theorem2}$ are similar.  The rest of this section is devoted to a proof of Proposition $\ref{higher regularity theorem}$, and we leave the proof of Proposition $\ref{higher regularity theorem2}$ as an exercise for the reader.

\begin{proof}
We start with the solution to the LANS equation $u$.  Then let $\delta>0$ be arbitrary, and let $v=t^\delta u$.  We note that $v(0)=0$.  Then 
\begin{equation*}\aligned \partial_t v&=\delta t^{\delta-1} u+t^\delta \partial_t u
\\ &=\delta t^{-1} v+t^\delta (\triangle u-\div(u\tensor u+\tau^\alpha (u,u)))
\\ &=\delta t^{-1} v+ \triangle v-t^{-\delta}\div (v\tensor v+\tau^\alpha (v,v)).
\endaligned
\end{equation*}
Applying Duhamel's principle, we get 
\begin{equation*}\aligned v&=e^{t\triangle}v_0+\int_0^t e^{(t-s)\triangle}s^{-1}v(s)ds+\int_0^t e^{(t-s)\triangle}s^{-\delta}V^\alpha(v(s),v(s))ds
\\ &=\int_0^t e^{(t-s)\triangle}v(s)ds+\int_0^t e^{(t-s)\triangle}V^\alpha(v(s),v(s))ds,
\endaligned
\end{equation*}
where the last line used that $v_0=0$.  Using $v=t^\delta u$, we get 
\begin{equation*}u=t^{-\delta}\int_0^t e^{(t-s)\triangle}s^{\delta-1}u(s)ds+t^{-\delta}\int_0^t e^{(t-s)\triangle}s^\delta V^\alpha(u(s),u(s))ds.
\end{equation*}

The key idea here is that we can choose $\delta$ to be large enough that both terms are integrable at $s=0$.  Now we are ready to apply the induction.  We have by assumption that the local solution $u$ is in $\dot{C}^T_{(s_2-s_1)/2;s_2,p}$, where $s_2\geq 1$.  For induction, we assume this solution $u$ is also in $\dot{C}^T_{(k-s_1)/2;k,p}$ where $k\geq s_2$, and seek to show that $u$ is in $\dot{C}^T_{(k+h-s_1)/2;k+h,\bar{p}}$, where $0<h<1$ is fixed and will be chosen later.  We have  
\begin{equation*}\aligned \|u\|_{H^{k+h,p}}&=t^{-\delta}\int_0^t \|e^{(t-s)\triangle} s^{\delta-1}u(s)\|_{H^{k+h,p}}ds + t^{-\delta}\int_0^t \|e^{(t-s)\delta} s^{\delta}V^\alpha(u(s))\|_{H^{k+h,p}}ds
\\ &\leq Ct^{-\delta}\int_0^t |t-s|^{-h/2}s^{\delta-1}\|u(s)\|_{H^{k,p}} + t^{\delta}\int_0^t |t-s|^{-(h+1+n/\tilde{p}-s_1)/2}s^{\delta}\|V^\alpha(u)\|_{H^{k-1}_{\tilde{p},q}} ds.
\endaligned
\end{equation*}

For the first piece, we have 
\begin{equation*}\aligned t^{-\delta}\int_0^t |t-s|^{-h/2}s^{\delta-1}\|u(s)\|_{H^{k,p}}&=t^{-\delta}\|u\|_{(k-s_1)/2;k,p}\int_0^t |t-s|^{-h/2}s^{\delta-1-(k-s_1)/2}ds
\\ &\leq C\|u\|_{(k-s_1)/2;k,p}t^{-\delta}t^{-h/2}t^{\delta-1-(k-s_1)/2+1}
\\ &\leq Ct^{-(k+h-s_1)/2}\|u\|_{(k-s_1)/2;k,p}
\endaligned
\end{equation*}
For this to work, we need the exponents of $|t-s|$ and $s$ in the integral to be strictly greater than negative $1$.  For $|t-s|$, this holds provided $h/2<1$.  For $s$, it works for a sufficiently large choice of $\delta$.  We note that without modifying the PDE to include these $t^\delta$ terms, we would need $(k-s_1)/2$ to be less than $1$, which does not hold for large $k$.

For the second piece, we split $V^\alpha$ into its two terms.  We begin by defining $\tilde{p}=p/2$, and we have 
\begin{equation*}\aligned &t^{-\delta}\int_0^t|t-s|^{-(h+n/\tilde{p}-s_1)/2}s^{\delta} \|\div(1-\triangle)^{-1}(\nabla u \nabla u)\|_{H^{k}_{\tilde{p},q}}ds 
\\ \leq &t^{-\delta}\int_0^t|t-s|^{-(h+s_1)/2}s^{\delta} \|(\nabla u \nabla u)\|_{H^{k-1}_{\tilde{p},q}}ds 
\\ \leq &t^{-\delta}\int_0^t|t-s|^{-(h+s_1)/2}s^{\delta} \|\nabla u \|_{H^{k-1,p}}\|\nabla u\|_{H^0_{p}}ds 
\\ \leq &t^{-\delta}\int_0^t|t-s|^{-(h+s_1)/2}s^{\delta} \| u \|_{H^{k,p}}\|u\|_{H^1_{p}}ds 
\\ \leq &t^{-\delta}\|u\|_{(k-s_1)/2;k,p}\|u\|_{(1-s_1)/2;1,p} \int_0^t|t-s|^{-(h+s_1)/2}s^{\delta-(k-s_1)/2-(1-s_1)/2}ds 
\\ \leq &t^{-(h+k-s_1)/2+1/2}\|u\|_{(k-s_1)/2;k,p}\|u\|_{H^{s_1,p}},
\endaligned
\end{equation*}
and this requires $\delta>(k-s_1)/2+(1-s_1)/2$ and $h+s_1<2$.  

For the last nonlinear piece, we have 
\begin{equation*}\aligned &t^{-\delta}\int_0^t|t-s|^{-(h+1+n/\tilde{p}-s_1)/2}s^{\delta} \|u\tensor u\|_{H^{k}_{\tilde{p},q}}ds 
\\ \leq &t^{-\delta}\int_0^t|t-s|^{-(h+1+s_1)/2}s^{\delta} \|u\|_{H^{k,p}}\|u\|_{H^0_{p}}ds 
\\ \leq&t^{-\delta}\|u\|_{(k-s_1)/2;k,p}\|u\|_{0;s_1,p} \int_0^t|t-s|^{-(h+1+s_1)/2}s^{\delta-(k-s_1)/2}ds 
\\ \leq &t^{-(h+k)/2+1/2}\|u\|_{(k-s_1)/2;k,p}\|u\|_{0;s_1,p}.
\endaligned
\end{equation*}
This requires $h+s_1<1$ (and thus $s_1<1$) and $\delta$ to be sufficiently large. 

So we finally get 
\begin{equation*}\aligned \|u\|_{H^{k+h,p}}&\leq C\|u\|_{(k-s_1)/2;k,p} (t^{-(k+h-s_1)/2}+t^{-(h+k-s_1)/2+(1-n/\bar{p})/2}+t^{-(h+k-s_1)/2+1/2})
\endaligned
\end{equation*}
which gives
\begin{equation*}\aligned &\|u\|_{(k+h-s_1)/2;k+h,p}
\\ \leq &C\|u\|_{(k-s_1)/2;k,p}t^{(k+h-s_1)/2}(t^{-(k+h-s_1)/2}+t^{-(h+k-s_1)/2+(1-n/\bar{p})/2}+t^{-(h+k-s_1)/2+1/2})
\\ \leq &C\|u\|_{(k-s_1)/2;k,p}(1+t^{(1-n/\bar{p})/2}+t^{1/2}).
\endaligned
\end{equation*}
This proves the desired result.  We remark that $\delta$ is chosen after beginning the induction step, while the appropriate value of $h$ is fixed by the choices of $n,$ $p,$ and $s_1$.
\end{proof}

\bibliographystyle{amsplain}
\bibliography{references2}

\end{document}